\documentclass[letterpaper, 10 pt, conference]{ieeeconf}  

\IEEEoverridecommandlockouts                              

\usepackage[english]{babel}
\usepackage{amsmath,amsfonts,stmaryrd,amssymb} 
\usepackage{bbding} 
\usepackage{enumerate} 
\usepackage{graphicx}
\usepackage{float} 
\usepackage{subcaption}
\usepackage{graphicx}
\usepackage{tabularx}

\usepackage[font=small]{caption}

\usepackage[english]{babel}

\usepackage{amsthm}

\newtheorem{definition}{Definition} 

\usepackage{geometry} 

\newtheorem{theorem}{Theorem}[section]

\newtheorem{remark}[theorem]{Remark}

\usepackage{algpseudocode}
\usepackage{algorithm}

\makeatletter
\newcounter{phase}[algorithm]
\newlength{\phaserulewidth}

\newcommand{\phase}[1]{%
  \vspace{-1.25ex}
  \Statex\leavevmode\llap{\rule{\dimexpr\labelwidth+\labelsep}{\phaserulewidth}}\rule{\linewidth}{\phaserulewidth}
  \Statex\strut\refstepcounter{phase}\textit{Phase~\thephase~--~#1}
  \vspace{-1.25ex}\Statex\leavevmode\llap{\rule{\dimexpr\labelwidth+\labelsep}{\phaserulewidth}}\rule{\linewidth}{\phaserulewidth}}
\makeatother

\usepackage{xcolor}
\usepackage{ulem}

\usepackage{hyperref}

\usepackage{graphicx}
\usepackage{wrapfig}

\usepackage{footmisc}

\newcommand{\mc}[1]{\mathcal{#1}}
\newcommand{\mbf}[1]{\mathbf{#1}}
\newcommand{\mbb}[1]{\mathbb{#1}}

\newcommand{\E}{\mathbb{E}}
\newcommand{\R}{\mathbb{R}}
\newcommand{\N}{\mathbb{N}}
\newcommand{\Z}{\mathbb{Z}}

\title{\LARGE \bf
Private and Accurate Decentralized Optimization via \\Encrypted and Structured Functional Perturbation
}

\author{Yijie Zhou and Shi Pu
\thanks{*This work was supported in parts by National Natural Science Foundation of China (NSFC) (Grant No. 62003287), Shenzhen Science and Technology Program (Grant No. RCYX20210609103229031 and No. GXWD20201231105722002-20200901175001001) and Shenzhen Institute of Artificial Intelligence and Robotics for Society (AIRS) (Grant No. AC01202101108).}
\thanks{Y. Zhou is with School of Data Science, The Chinese University of Hong Kong, Shenzhen, China. S. Pu is with School of Data Science, Shenzhen Institute of Artificial Intelligence and Robotics for Society (AIRS), The Chinese University of Hong Kong, Shenzhen, China 
        ({\tt\small zyjbuaa@gmail.com, pushi@cuhk.edu.cn}).}%
}

\begin{document}

\maketitle
\thispagestyle{empty}
\pagestyle{empty}

\def\sp#1{\color{blue}{#1}}
\def\alert#1{\color{red}{#1}}

\begin{abstract}
We propose a decentralized optimization algorithm that preserves the privacy of agents' cost functions without sacrificing accuracy, termed EFPSN. 
The algorithm adopts Paillier cryptosystem to construct zero-sum functional perturbations. 
Then, based on the perturbed cost functions, any existing decentralized optimization algorithm can be utilized
to obtain the accurate solution.
We theoretically prove that EFPSN is ($\epsilon, \delta$)-differentially private and can achieve nearly perfect privacy under deliberate parameter settings. Numerical experiments further confirm the effectiveness of the algorithm.
\end{abstract}

\section{Introduction}

The problem of optimizing a global objective function through the cooperation of multiple agents has gained increased attention in recent years \cite{nedic2009distributed,pu2021distributed}. This is driven by the wide applicability of the problem to many engineering and scientific domains, ranging from cooperative control, distributed sensing, multi-agent systems, and sensor networks to large-scale machine learning; see, e.g., \cite{pu2016noise,yan2008design, murphey2002cooperative, roy2019braintorrent}.

In this paper, we consider a peer-to-peer network of $n$ agents that solves the following optimization problem cooperatively:
\begin{equation}
    \label{eqn:problem_definition}
    \min_{x \in \mc{D}} F(x) = \frac{1}{n} \sum_{i\in \mc{N}} f_i(x),
\end{equation}
where $x$ is the common parameter of all the agents, $\mc{D}$ is the domain of $x$, and $\mc{N}$ denotes the collection of all the agents.

Despite the enormous success of gradient-based distributed optimization algorithms, they all require agents to explicitly share optimization variables and/or gradient estimates in every iteration. This would become a problem in applications involving sensitive data. Zhu et al. \cite{zhu_deep_nodate} showed that obtaining private training data from publicly shared gradients is possible, and the recovery is pixelwise accurate for images and token-wise matching for texts. Consequently, we wish to solve (\ref{eqn:problem_definition}) in a private way.

In the machine learning community, one of the common forms of the individual function is $f_i(x)\triangleq \E_{\xi_i \sim \mc{D}_i} [l_i (x, \xi_i)]$, where $\mc{D}_i$ is the local data distribution of agent $i$, and $\xi_i$ is a data sample or a batch of data samples. In such a case, each function $f_i$ contains information on the data distribution $\mc{D}_i$, which is usually sensitive. It is thus crucial to keep the objective functions private.
More specifically, by privacy, we refer to keeping all $f_i$ {from being inferred by} adversaries. In this work, we consider two types of adversaries:

\begin{itemize}
    \item The eavesdropper: an external adversary having access to all information transmitted through the communication channels within the network.
    \item The curious-but-honest adversary: an external adversary that corrupts a subset of the agents. The adversary knows all the information of every corrupted agent $i$, including all the information within agent $i$, e.g., the individual function $f_i$, and all the information passed from its neighbor to agent $i$. However, each corrupted agent obeys the optimization protocol precisely.
\end{itemize}

Plenty of efforts have been reported to counteract potential privacy breaches in distributed optimization. Privacy preserving algorithms either process the messages being transmitted between agents or change the functions to be optimized. We refer to them as message-based and function-based methods,  respectively. Differential Privacy (DP)~\cite{cortes_differential_2016, dwork_algorithmic_2013}, a de facto standard for privacy serving, has been introduced to the context of distributed optimization~\cite{chen_differential_2021, lou_privacy_2015, huang_differentially_2015, nozari_differentially_2016}. The DP-based methods can be categorized into message-perturbing \cite{chen_differential_2021, lou_privacy_2015, huang_differentially_2015} and function-perturbing ones \cite{nozari_differentially_2016}. The message-perturbing methods inject noises to the messages each agent sends, while the latter adds functional noises into each agent's cost function $f_i$. 

However, the direct combination of DP and distributed optimization suffers from the accuracy-privacy trade-off. Huang et al. \cite{huang_differentially_2015} observed that, by fixing the other parameters, the accuracy level of the obtained solution is in the order of $O (\frac{1}{\epsilon^2})$, where $\epsilon$ is the privacy budget that is inversely proportional to privacy. The papers \cite{lu_privacy_2018, zhang_enabling_2019} adopted encryption techniques to preserve privacy in distributed optimization. Nevertheless, the heavy communication and computation overhead prevents the methods from many real-world applications. 

Adding structured noise \cite{gupta_preserving_2021} is a workaround for the accuracy-privacy trade-off and the heavy overhead mentioned above. The paper \cite{gupta_preserving_2021} constructed a set of zero-sum Gaussian noises to perturb the affine terms of the objective functions that are assumed to be quadratic. Nonetheless, the method fails under an eavesdropping attack due to the naive communication method. Besides, the privacy analysis in \cite{gupta_preserving_2021} is carried out under a self-defined privacy framework.
We provide a categorization of the current privacy-preserving distributed optimization algorithms in Table \ref{tab:all_alg}.

\begin{table}[]
    \centering
    \scalebox{0.85}{
    \begin{tabular}{||c|c|c|c|c|c||}
    \hline
       Paper Index  & \cite{chen_differential_2021, lou_privacy_2015,huang_differentially_2015} & \cite{nozari_differentially_2016} & \cite{lu_privacy_2018,zhang_enabling_2019} & \cite{gupta_preserving_2021} & Ours\\
        \hline
        Message-based &  \Checkmark & \XSolidBrush & \Checkmark  & \XSolidBrush & \XSolidBrush\\
        \hline
        Function-based & \XSolidBrush & \Checkmark & \XSolidBrush & \Checkmark & \Checkmark \\
        \hline
   Structured-noise & \XSolidBrush & \XSolidBrush & \XSolidBrush & \Checkmark & \Checkmark \\
        \hline
        Encryption & \XSolidBrush & \XSolidBrush & \Checkmark & \XSolidBrush & \Checkmark\\
        \hline
              DP-based & \Checkmark & \Checkmark & \XSolidBrush & \XSolidBrush & \Checkmark\\
        \hline
        
    \end{tabular}}
    \caption{\label{tab:all_alg} Categorization of existing privacy-preserving distributed optimization algorithms.}
    \label{tab:my_label}
\end{table}

In this paper, we integrate the encryption-based scheme and the structured noise method deliberately under the functional DP framework. The proposed new method, termed the Encrypted Functional Perturbation with Structured Noise algorithm (EFPSN), enjoys the benefits of several previous methods. In particular, EFPSN adopts the Paillier encryption scheme to construct zero-sum noises among the agents secretly. The noises are subsequently used to generate functional perturbations for each agent's cost function. Such a
procedure differs from those in \cite{lu_privacy_2018, zhang_enabling_2019} and bypasses the heavy communication and computation overhead caused by encryption at every iteration. Then, based on the perturbed cost functions, any existing decentralized optimization algorithm can be utilized
to obtain the accurate solution to problem \eqref{eqn:problem_definition} thanks to the structured noises. We further theoretically prove that EFPSN is ($\epsilon, \delta$)-differentially private and can achieve nearly perfect privacy under deliberate parameter settings.
In other words, EFPSN can achieve nearly perfect privacy without sacrificing accuracy.

The rest of this paper is organized as follows: Section 2 specifies the notation and provides related background knowledge. Then, we develop the Encrypted Functional Perturbation with Structured Noise algorithm in Section 3. Privacy analysis under the DP framework is carried out in Section 4. Finally, simulation examples and conclusion are presented in Section 5 and Section 6 correspondingly.


\section{Preliminaries}

This section introduces the notation, graph-related concepts and some background knowledge on Hilbert space and Paillier Cryptosystem, since generating functional perturbations relies on the orthonormal systems in some Hilbert space, and Paillier Cryptosystem \cite{paillier1999public} is adopted to construct structured noises privately.

\subsection{Notation}

We use $\R$ to denote the set of real numbers and $\R^d$ the Euclidean space of dimension $d$. The space of scalar-valued inﬁnite sequences are denoted by $\R^\N$. Let $\Z,\Z_{>0}$ be the set of integers and positive integers, respectively. Given $w\in \Z_{>0}$ and $\Z_w \triangleq \{0,1,\cdots,w\}$, $\Z_w^\ast$ denotes the set of positive integers which are smaller than $w$ and do not have common factors other than 1 with $w$.
Let $l_2\subset \R^{\N}$ be the space of infinite square summable sequences. For $D\subseteq \R^d$, $L_2(D)$  denotes the set of square-integrable measurable functions over $D$.
$\mbf{1}$ denotes a column vector with all entries equal to 1. A vector is viewed as a column vector unless otherwise stated. $A^T$ denotes the transpose of the matrix $A$ and $x^Ty$ denotes the scalar product of two vectors $x$ and $y$. We use $\langle \cdot, \cdot \rangle$ to denote the inner product and $||\cdot||$ to denote the Euclidean norm for a vector (induced Euclidean norm for a matrix). A square matrix $A$ is column-stochastic when its elements in every column add up to one. A matrix $A$ is said to be doubly stochastic when both $A$ and $A^T$ are column-stochastic matrices. 

We use $\mbb{P}\{\mc{A}\}$ to denote the probability of an event $\mc{A}$, $\mc{P}_{X}(y)$ the probability density function of random variable $X$, and $\E [X | \mc{F}]$ the expectation of a random variable $X$ with $\mc{F}$ denoting the sigma algebra, which will be omitted when clear from the context. For an encryption scheme, $\text{En}(\cdot), \text{De}(\cdot)$ represent the encoder and decoder respectively. Let $\gcd, \text{lcm}, \mod$ be the greatest common divider, the least common multiple, and the modulo operator, respectively. $N(\mu,\sigma^2)$ is the (multivariate) Gaussian distribution with (vector) mean $\mu$ and variance (covariance matrix) $\sigma^2$. $N^\dagger$ represents the degenerate Gaussian distribution. 

\subsection{ Graph Related Concepts}

We assume that the agents interact on an undirected graph, described by a matrix $W\in\mathbb{R}^{n\times n}$. More specifically, if agents $i$ and $j$ can communicate and interact with each other, then $w_{ij}$, the $(i,j)$-th entry of $W$, is positive. Otherwise, $w_{ij}$ equals zero. The neighbor set $\mc{N}_i$ of agent $i$ is defined as the set of agents $\{j|w_{ij}>0\}$. Note that $i\in \mc{N}_i$ always holds. Denote $\mc{L}$ as the graph Laplacian depicted by $W$. Let $\mu_1 \le \mu_2 \le \cdots \le \mu_n$ be the eigenvalues of $\mc{L}$ and $M$ be the unitary matrix that satisfies $\mc{L} = M \text{Diag} (\mu_1,..., \mu_n) M^T$. Denote $\underline{\mu}(\mc{L})$ as the second smallest eigenvalue of $\mc{L}$ and $\bar{\mu}(\mc{L})$ as the largest eigenvalue of $\mc{L}$. 

\subsection{Hilbert Spaces}

A Hilbert space $\mc{H}$ is a complete inner-product space. A set $\{e_k\}_{k\in \N} \subset \mc{H}$ is an orthonormal system if $\langle e_k, e_j\rangle = 0$ for $k\ne j$ and $\langle e_k, e_k\rangle = ||e_k||^2=1$ for all $k\in \N$. If, in addition, the set of linear combinations of $\{e_k\}_{k \in \N}$ is dense in $\mc{H}$, then $\{e_k\}_{k\in I}$ is an orthonormal basis. If $\mc{H}$ is separable, then any orthonormal basis is countable, and we have
\begin{equation}
    \label{eqn:hilbert_space}
    h = \sum_{k=0}^\infty \langle h, e_k \rangle e_k,
\end{equation}
for any $h\in \mc{H}$. Define the coefficient sequence $\mbb{\theta} \in \R^{\N}$ by $\theta_k = \langle h, e_k \rangle$ for $k\in \N$. Then $\mbb{\theta} \in l_2$ and, by Parseval's identity, $||h|| = ||\mbb{\theta}||$. Let $\Phi: l_2 \rightarrow \mc{H}$ be the linear bijection that maps the coefficient sequence $\mbb{\theta}$ to $h$. For an arbitrary $D\subseteq \R^d, L_2(D)$ is a Hilbert space, and the inner product is the integral of the product of functions. Moreover, $L_2(D)$ is separable. In this paper, we denote $\{e_k\}_{k=0}^\infty$ as an orthonormal basis for $L_2(D)$ and $\Phi:l_2 \rightarrow L_2(D)$ the corresponding linear bijection between coefficient sequences and functions.

\subsection{Paillier Cryptosystem}

The Paillier Cryptosystem is an algorithm for public key cryptography. The algorithm applies to the scenario of sending a private message over open and insecure communication links, and it consists of key generation, encryption, and decryption steps as follows:

\begin{itemize}
    \item Key Generation
    \begin{itemize}
        \item The message receiver chooses two large prime numbers $a$ and $b$ randomly and independently of each other such that $\gcd(ab,(a-1)(b-1))=1$. This property is assured if both primes are of equal length \cite{fulton_2010}.
        \item Compute $f=ab$ and $\lambda = \text{lcm}(a-1, b-1)$
        \item Seletect random integer $g$, where $g\in \Z_{f^2}^\ast$ such that the modular multiplicative inverse $\mu = (\frac{(g^\lambda \mod f^2)-1}{f})^{-1} \mod f$ exists.
        \item Let the public key $\bar{\mc{K}}=(f,g)$, the private key $\tilde{\mc{K}}=(\lambda, \mu)$.
    \end{itemize}
    \item Encryption: To encrypt a plaintext $\underline{p} \in \Z_f$, the sender selects a random number $r \in \Z_f ^\ast$ and computes the ciphertext $\underline{c} = \text{En} (\underline{p}, \bar{\mc{K}}, r) = g^{\underline{p}} \cdot r^f \mod f^2$.
    
    \item Decryption: To decrypt the ciphertext $\underline{c} \in \Z_{f^2}$, the receiver computes the decrypted text $\bar{\underline{p}} = \text{De} (\underline{c}, \bar{\mc{K}}, \tilde{\mc{K}}) = \frac{(\underline{c}^{\lambda} \mod f^2)-1}{f} \cdot \mu \mod f$.
\end{itemize}

One notable homomorphic property the Paillier encryption scheme holds is: given any $\underline{p}_1, \cdots, \underline{p}_m \in \N$. If $\sum_{l=1}^m \underline{p}_l \in Z_f, \text{then } \text{De}(\prod_{l=1}^m \text{En}(\underline{p}_l, \bar{\mc{K}}, r_l), \bar{\mc{K}}, \tilde{\mc{K}}) = \sum_{l=1}^m \underline{p}_l$.


\section{Algorithm Design}

In this section, we propose the Encrypted Functional Perturbation with Structured Noise algorithm (EFPSN in short) that privately solves the problem (\ref{eqn:problem_definition}). Unlike the majority of privacy-preserving algorithms, EFPSN does not sacrifice accuracy for privacy.

To achieve privacy, EFPSN adds structured functional perturbations on individual cost functions $\{f_i\}_{i\in N}$. Specifically, the algorithm aims at making the functional perturbations zero-sum. EFPSN consists of two phases, as shown in Algorithm \ref{alg:EFPSN}. 

\begin{algorithm}
\caption{Encrypted Functional Perturbation with Structured Noise}\label{alg:EFPSN}
\begin{algorithmic}[1]
\Require Cost function $\{f_i\}_{i\in \mc{N}}$, noise precision order $P$, perturbation order $K$, and $\{\sigma_k\}_{k\in \Z_K}$.
\Ensure $x^*$

\phase{Masking cost functions}
\For{$i\in\mc{N}$}
    \State {Generate key pair $(\bar{\mc{K}}_i, \tilde{\mc{K}}_i)$, and $r_i$ following the Paillier encryption scheme}
    \State {Share public key $\bar{\mc{K}}_i$ with agent $j\in \mc{N}_i$}
\EndFor
\For{$i\in\mc{N}$}
    \For{$(j,k)\in\mc{N}_i \times \Z_K$}

        \State{Generate random noise $\eta_{ijk} \sim N(0,\sigma_k^2)$}
        
        \State{Calculate $\underline{c}_{ijk}= \text{En}(\lfloor 10^P \eta_{ijk}\rfloor,\bar{\mc{K}}_j,r_i)$}
        
        \State{Send $\underline{c}_{ijk}$ to agent $j$}
        
    \EndFor
    
    \For{$k\in\Z_K$}
        \State {Calculate $\bar \eta_{ik} = \sum_{j\in\mc{N}_i} \eta_{ijk} - $
        
        $\quad\quad10^{-P}\text{De}(\prod_{j\in\mc{N}_i} \underline{c}_{jik}, \bar{\mc{K}}_i, \tilde{\mc{K}}_i)$ }
    \EndFor

\State {$\hat f_i = \Phi (\Phi^{-1} (f_i) + \bar \eta_i)= f_i + \Phi (\bar \eta_i), $

$\quad \bar \eta_i = [\bar \eta_{i0}, ..., \bar \eta_{iK},0,...]\in\R^{\N}$ }

\EndFor

\phase{Distributed optimization}
\State{Execute any distributed optimization algorithm on the masked functions $\{\hat{f}_i\}_{i\in \mc{N}}$}

\end{algorithmic}
\end{algorithm}

In phase I, the agents in the network cooperate to generate functional perturbations in a way that is immune to eavesdropping attacks and honest-but-curious attacks to some extent. First, they generate keys and random numbers required by the Paillier encryption scheme. Then, they encrypt and send random noise $\lfloor 10^P \eta_{ijk}\rfloor$ to their neighbors, and the noises are further decrypted to construct the zero-sum perturbation. Due to encryption, the signals are transferred privately and securely under eavesdropping. However, since Paillier encryption only works for integers, we set a precision order $P$, and encrypt and decrypt $\lfloor 10^P \eta_{ijk}\rfloor$.

Subsequently, in Line 10, each agent $i$ calculates $\bar{\eta}_{ik}$ by subtracting the sum of noises it receives from the sum of noises it sends. Due to Paillier encryption's homomorphic property, each agent only needs to decode once for each $k \in \Z_K$. This saves computation, especially when each agent has a large number of neighbors.

Paillier encryption scheme guarantees privacy under eavesdropping attacks. In terms of the honest-but-curious attacks, the noise coefficient sequence $\bar\eta_i$ of each agent $i\in\mc{N}$ remains unknown to the attacker as long as the attacker does not corrupt all $i$'s neighbors. Under such circumstances, the privacy of agent $i$ is maintained.

It is worth noting that $\sum_{i\in \mc{N}} \sum_{j\in \mc{N}_i} (\eta_{ijk} - \eta_{jik}) = 0$ for all $k$. Such a construction forces $\lim_{P \rightarrow \infty} \sum_{i} \bar{\eta}_{ik} = 0$ to hold for all $k$. Therefore, we have generated a set of zero-sum signals $\{\bar{\eta}_{ik}\}_{i\in \mc{N}}$ given large $P$.

Finally in Line 12, the agents perturb the cost functions by adding $\Phi(\bar\eta_i)$. $\Phi(\cdot)$ is a function that maps a sequence in $l_2$ to a function in  $L_2(D)$. Such a construction depends on the orthonormal system $\{e_k\}_{k\in I}$ in $\mc{H}$ one chooses. For instance, given the orthonormal system $\{e_k\}_{k\in I}$ and a sequence $\bar\eta_i \triangleq \{\bar\eta_{ik}\}_{k \in I}$, $\Phi(\bar\eta_i) = \sum_{k\in I} \bar\eta_{ik} e_k$. The orthonormal system we use will be specified later. 

Since $\{\bar{\eta}_{i}\}_{i\in \mc{N}}$ is zero-sum when $P \rightarrow \infty$, we have
\begin{equation}
    \begin{aligned}[b]
    \label{eqn:accuracy}
    \lim_{P\rightarrow \infty}\sum_i \Phi(\bar\eta_i) 
    &= \lim_{P\rightarrow \infty} \sum_i \sum_k \bar\eta_{ik} e_k\\
    &= \sum_k \lim_{P\rightarrow \infty} (\sum_i \bar\eta_{ik}) e_k\\
    &= 0.
    \end{aligned}
\end{equation}
Therefore, the set of perturbing functions $\{\Phi(\bar\eta_i)\}_{i\in \mc{N}}$ is zero-sum when $P \rightarrow \infty$. Additionally, the decrypted text in Line 10 is of precision $10^{-P}$. Consequently, the error brought by $P$ will be dominated by the floating-point error once $P$ is set to a moderately large value. Though $\{\Phi(\bar\eta_i)\}_{i\in \mc{N}}$ is zero-sum, each agent $i$ gains privacy from the non-zero functional perturbation $\Phi(\bar\eta_i)$.

In phase II, the agents may conduct any distributed optimization algorithm on $\{\hat f_i\}_{i\in\mc{N}}$. Since $\sum_i \hat f_i(x) = \sum_i f_i(x)$, the obtained solution solves the original problem (\ref{eqn:problem_definition}) when $P \rightarrow \infty$. Namely, EFPSN solves problem (\ref{eqn:problem_definition}) without any accuracy degeneration given proper $P$.

\begin{remark}
EFPSN combines encryption, functional perturbation, and structured noise and is superior to considering any one of the techniques. In previous message-based methods, encryption at every iteration results in heavy communication and computation overhead, whereas the function-level encryption in EFPSN alleviates such a pain: only insignificant communication and computation overhead occur in phase I. Regarding the existing function-based methods, the obtained solution after functional perturbation suffers from a privacy-related deviation from the solution to problem \eqref{eqn:problem_definition}, and thus leads to privacy-accuracy tradeoff. For EFPSN, however, the optimization accuracy is independent of the privacy budget, which will be elaborated more in the following sections.
Moreover, using structured noise alone fails with the presence of eavesdropping, limiting its applicability.
\end{remark}


\section{Privacy Analysis}
In this section, we analyse the privacy-related property of EFPSN under the framework of differential privacy.
Particularly, we prove the mechanism of generating the masked functions is differentially private. And since DP is compatible with any post-processing, EFPSN remains differentially private. 

First, we introduce the definition of $\mc{V}$-adjacency, which was originally proposed in \cite{nozari_differentially_2016}.

\begin{definition} [$\mc{V}$-adjacency]
\label{def:adjacency}
Given any normed vector space $(\mc{V}, ||\cdot||_{\mc{V}})$ with $\mc{V} \subseteq L_2(D)$, two sets of functions $F=\{f_i\}_{i\in\mc{N}},F'=\{f'_i\}_{i\in\mc{N}} \subset L_2(D)$ are $\mc{V}$-adjacent if there exists $I\in \mc{N}$ such that 

\begin{equation}
    f_i = f_i', i\ne I, \text{ and } f_{I} - f_{I}' \in \mc{V}.
\end{equation}
\end{definition}

We adopt the standard $(\epsilon,\delta)$-DP definition under our functional setting.

\begin{definition} [$(\epsilon,\delta)$-Differential Privacy]
\label{def:privacy}
Consider a random map $\mc{M}: L_2(D)^n\rightarrow \mc{X}$ from the function space $L_2(D)^n $ to an arbitrary set $\mc{X}$. Given $\epsilon,\delta \in \mbb{R}_{\ge 0}$, the map $\mc{M}$ is $(\epsilon,\delta)$-differentially private if, for any two $\mc{V}$-adjacent sets of functions $F=\{f_i\}_{i\in\mc{N}}$ and $F'=\{f'_i\}_{i\in\mc{N}}$, one has
\begin{equation}
    \mbb{P} \{\mc{M} (F) \in \mc{O} \} \le e^{\epsilon} \cdot \mbb{P} \{\mc{M} (F') \in \mc{O} \} + \delta.
\end{equation}
\end{definition}

We choose our adjacency space $\mc{V}_q$ as follows. Given $q>1$, consider the weight sequence $\{k^q\}^\infty_{k=1}$ and define the adjacency vector space to be the image of the resulting weighted $l_2$ space under $\Phi$, i.e.,

\begin{equation}
\begin{aligned}[b]
    \mc{V}_q = \Phi(\{\delta\in \mbb{R}^{\mbb{N}} | \sum_{k=1}^\infty k^{2q} \delta_k^4<\infty\}),
\end{aligned}
\label{eqn:adj_space}
\end{equation}
where $\delta_k$ is the $k$th element  of $\delta$. The rationale of considering such a space will be made clear from the analysis. Moreover, 
\begin{equation}
\begin{aligned}[b]
    ||f||_{\mc{V}_q}\triangleq (\sum_{k=1}^\infty(k^{2q}\delta_k^4))^{\frac{1}{4}}, \text{with } \delta = \Phi^{-1}(f)
\end{aligned}
\label{eqn:adj_norm}
\end{equation}
is a norm on $\mc{V}_q$.

Now we introduce our main theorem about the privacy-preserving property of EFPSN.

\begin{theorem}\label{theo:privacy}
Given $q>1,\gamma>0, p\in(1/2,q-1/2)$, the chosen $\mc{V}_q$ space, and $\sigma_k^2 =  \frac{\gamma}{k^p}$, the mechanism in Alg. \ref{alg:EFPSN} is $(\epsilon, \delta)$-differential private when the precision order and the perturbation order $P, K \rightarrow \infty$, with $\epsilon = \left\{ \frac{1}{\underline\mu(\mc{L})}(\frac{A}{4} + \frac{R\sqrt{\bar{\mu}(\mc{L})A}}{\sqrt{2}})   \right\}, \delta=e^{-\frac{R^2}{2}}$, where $A=\frac{1}{\gamma} \sqrt{\zeta (2(q-p))} ||f_I - f_I'||^2_{\mc{V}_q}$, and $R$ is an arbitrary positive real number.
\end{theorem}

\begin{proof}
When the precision order $P \rightarrow \infty$, the encryption process is precise. For convenience, we can ignore the encryption process and the flooring in Alg. \ref{alg:EFPSN}. Denote
\begin{equation}
    \begin{aligned}
        \eta_{ik} &\triangleq \sum_{j\in\mc{N}_i} \eta_{ijk} - \sum_{j\in\mc{N}_i} \eta_{jik},\\
        \eta_i &\triangleq [\eta_{i0},...,\eta_{iK}].
    \end{aligned}
\end{equation}
Under such a case, $\mc{M}$ in algorithm \ref{alg:EFPSN} is essentially adding functional perturbations constructed from a set of degenerate Gaussian noises. Specifically, given $F=\{f_i\}_{i\in \mc{N}}$, we have
\begin{equation}\label{eqn:mech}
    \mc{M}(F) = \{f_i + \Phi(\eta_i)\}_{i \in \mc{N}},
\end{equation}
and 
\begin{equation}\label{eqn:Gau}
    \pmb{\eta}_k=[\eta_{1k},...,\eta_{nk}]^T \sim N^\dagger(0_n,2 \sigma_k ^2 \mc{L}), \forall k \in \Z_{K}
\end{equation}
are the zero-sum coefficients. 

From \cite{gupta_preserving_2021}, we have 
\begin{equation}
\mc{P}_{\pmb{\eta}_k}(\pmb y) = \left\{
\begin{aligned}
\frac{1}{\sqrt{\det ^*(4\pi \sigma_k^2 \mc{L})}}\exp (-\frac{\pmb y^T\mc{L}^\dagger \pmb y}{4\sigma_k^2})&, \pmb y^T\pmb{1}=0 \\
0&,\text{else},
\end{aligned}
\right.
\label{eqn:deg_gau}
\end{equation}
where $\mc{L}^\dagger = M \text{Diag}(0, 1/\mu_2,...,1/\mu_n) M^T$ and $\det^*(4\pi\sigma_k^2\mc{L}) = (4\pi \sigma_k^2)^{n-1}\prod_{i=2}^n \mu_i$.

Let $ F'=\{f_i'\}_{i\in \mc{N}}$ be a set of $\mc{V}$-adjacent functions of $F$ that differ only in the $I$-th element. 
 Let $\Psi^{-1}:L_2(D)^n \rightarrow \R^{n\times \N}$ be the map such that $\Psi^{-1}(F)=\{\Phi^{-1}(f_i)\}_{i\in\mc{N}}$.
Define $\Phi_{0:k}^{-1}: L_2(D)\rightarrow \mbb{R}^{k+1}$, $\Phi_{k}^{-1}: L_2(D)\rightarrow \mbb{R}$ as the map that returns the first $k+1$ coefficients and the $
(k+1)$-th coefficient of $\Phi^{-1}(\cdot)$, respectively. Similarly, define $\Psi_{0:k}^{-1}: L_2(D)^n \rightarrow \R^{n\times (k+1)}, \Psi_{k}^{-1}: L_2(D)^n \rightarrow \R^{n}$ as the map that returns the first $k+1$ columns and the $
(k+1)$-th column of $\Psi^{-1}(\cdot)$, respectively. For any $\mc{O}\in L_2(D)^n$, $\mc{O}_i$ is the $i$th element of $\mc{O}$. Let $\mc{O}_i - f_i \triangleq \{g_i \in L_2(D)| g_i + f_i \in \mc{O}_i\}$ and $\mc{O} - F \triangleq \{\{g_i\}_{i\in\mc{N}} \in L_2(D)^n| \{g_i + f_i\}_{i\in \mc{N}} \in \mc{O}\}$.

We have
\begin{equation}
\begin{aligned}[b]
    &\mbb{P}\{\mc{M}(F)\in\mc{O}\} \\
    &= \mbb{P}\{\{\eta_i\}_{i\in\mc{N}}\in \Psi^{-1}(\mc{O}- F)\}\\
    &= \lim_{K\rightarrow \infty} \int_{\Psi_{0:K}^{-1}(\mc{O} - F)} \prod_{k=0}^K \mc{P}_{\pmb\eta_k}(\pmb y_k) d\pmb y_0\cdots d\pmb y_K
\end{aligned}
\label{eqn:prob:M(F)}
\end{equation}
and 
\begin{equation}
\begin{aligned}[b]
    &\mbb{P}\{\mc{M}(F')\in\mc{O}\} \\
    &= \mbb{P}\{\{\eta_i\}_{i\in\mc{N}}\in \Psi^{-1}(\mc{O}- F')\}\\
    &= \lim_{K\rightarrow \infty} \int_{\Psi_{0:K}^{-1}(\mc{O} - F')} \prod_{k=0}^K \mc{P}_{\pmb\eta_k}(\pmb y_k) d\pmb y_0\cdots d\pmb y_K
\end{aligned}
\label{eqn:prob:M(F')}
\end{equation}

 By the linearity of $\Phi$, we have $\Phi^{-1}(\mc{O}_i - f'_i) = \Phi^{-1} (\mc{O}_i - f_i) + \xi_i $ for all $\mc{O} \subseteq L_2(D)^n$ and $i\in\mc{N}$, where $\xi_i \triangleq \Phi^{-1}( f_i - f'_i) $. Denoting $\pmb\xi_k=[{\xi}_{1k},...,{\xi}_{nk}]^T$, we have 

\begin{equation}
    \Psi_k^{-1}(\mc{O} - F') = \Psi_k^{-1}(\mc{O} - F) + \pmb\xi_k
\label{eqn:prob:diff}
\end{equation}

Combing (\ref{eqn:prob:M(F')}) and (\ref{eqn:prob:diff}), we have

\begin{equation}
\begin{aligned}[b]
    &\mbb{P}\{\mc{M}(F')\in\mc{O}\}\\
    &=\lim_{K\rightarrow \infty} \int_{\Psi_{0:K}^{-1}(\mc{O} - F)} \prod_{k=0}^K \mc{P}_{\pmb\eta_k}(\pmb y_k+\pmb\xi_k) d\pmb y_0\cdots d\pmb y_K
\end{aligned}
\label{eqn:prob:M(F')new}
\end{equation}

To prove that $\mc{M}$ is $(\epsilon, \delta)$-DP, it suffices to show the ratio of $\prod_{k=0}^K \mc{P}_{\pmb\eta_k}(\pmb\eta_k)$ over $\prod_{k=0}^K \mc{P}_{\pmb\eta_k}(\pmb\eta_k+\pmb\xi_k)$ is bounded by $e^{\epsilon}$ with probability at least $(1-\delta)$.

We know that
\begin{equation}
\begin{aligned}[b]
    &\prod_{k=0}^K\frac{\mc{P}_{\pmb\eta_k}(\pmb\eta_k)}{\mc{P}_{\pmb\eta_k}(\pmb\eta_k+\pmb\xi_k)} \\
    &= \exp\left\{ \sum_{k=0}^K \frac{2\pmb\xi_k^T \mc{L}^\dagger \pmb\eta_k + \pmb\xi_k^T \mc{L}^\dagger \pmb\xi_k}{4\sigma_k^2} 
    \right\}\\
    &\le \exp\left\{\frac{1}{\underline \mu (\mathcal{L})} (\sum_{k=0}^K \frac{\pmb\xi_k^T   \pmb\eta_k}{2\sigma_k^2} + \sum_{k=0}^K\frac{||\pmb\xi_k||^2}{4\sigma_k^2})\right\}\\ 
\end{aligned}
\label{eqn:prob:ratio}
\end{equation}


Denote
$$\textbf{Rat}\triangleq  \exp\left\{\frac{1}{\underline \mu (\mathcal{L})} (\sum_{k=0}^\infty \frac{\pmb\xi_k^T   \pmb\eta_k}{2\sigma_k^2} + \sum_{k=0}^\infty\frac{||\pmb\xi_k||^2}{4\sigma_k^2})\right\}.$$
We show that $\textbf{Rat}$ is bounded with certain probability.

Since $\{f\}_{i\in\mc{N}}$ and $\{f'\}_{i\in\mc{N}}$ only differ in one element, $\pmb\xi_{k}$ has at most one non-zero element. Noting that $\pmb\eta_k$ is random and drawn from $N^\dagger(0_n,2 \sigma_k ^2 \mc{L})$, $\frac{\pmb\xi_k^T  \pmb\eta_k}{2\sigma_k^2}$ is a univariate Gaussian random variable.  From (\ref{eqn:Gau}), each element of $\pmb\eta_k$ is at most of variance $2\sigma_k^2 \bar \mu_\mc{L}$. Thus, $\sum_{k=0}^\infty \frac{\pmb\xi_k^T  \pmb\eta_k}{2\sigma_k^2}$ is a univariate Gaussian random variable with variance less than or equal to $\frac{\bar \mu(\mc{L})}{2}\sum_{k=0}^\infty\frac{||\pmb\xi_{k}||^2 }{\sigma^2_k}$. 

We further bound the summation term in its variance following a similar way in the proof of Theorem V.2 in \cite{nozari_differentially_2016}:
\begin{equation}
\begin{aligned}[b]
    \sum_{k=0}^\infty \frac{||\pmb\xi_k||^2}{\sigma_k^2} 
    &= \sum_{k=0}^\infty
        \frac{k^q||\pmb \xi_k||^2}{k^q\sigma_k^2} \\
    &\le (\sum_{k=0}^\infty \frac{1}{(k^q \sigma_k^2)^2})^{\frac{1}{2}} (\sum_{k=0}^\infty (k^q ||\pmb\xi_k||^2)^2)^{\frac{1}{2}} \\
    & = \frac{1}{\gamma} \sqrt{\zeta (2(q-p))} ||f_I - f_I'||^2_{\mc{V}_q}\\
    &\triangleq A.
\end{aligned}
\label{eqn:bound_noise}
\end{equation}

Let $R$ be an arbitrary positive real number. When $\sum_{k=0}^\infty \frac{\pmb\xi_k^T  \pmb{\eta}_k}{2\sigma_k^2} \le  R\sqrt{\frac{\bar \mu(\mc{L})}{2}\sum_{k=0}^\infty\frac{||\pmb\xi_{k}||^2 }{\sigma^2_k}}$ holds, we have

\begin{equation}
\begin{aligned}[b]
   \textbf{Rat}
   &\le \exp \left\{\frac{1}{\underline\mu(\mc{L})} (R \sqrt{\sum_{k=0}^\infty\frac{||\pmb\xi_{k}||^2 }{\sigma^2_k}} + \frac{1}{4} \sum_{k=0}^\infty\frac{||\pmb\xi_{k}||^2 }{\sigma^2_k})\right\}\\
   &\le \exp\left\{ \frac{1}{\underline\mu(\mc{L})}(\frac{A}{4} + \frac{R\sqrt{\bar{\mu}(\mc{L})A}}{\sqrt{2}})   \right\}.
\end{aligned}
\label{eqn:ratio_bound}
\end{equation}

By adopting the Chernoff bound for Gaussian random variable, we have

\begin{equation}
\begin{aligned}[b]
    \mbb{P}\left\{\sum_{k=0}^\infty \frac{\pmb\xi_k^T  \pmb \eta_k}{2\sigma_k^2} \ge   R\sqrt{\frac{\bar \mu(\mc{L})}{2}\sum_{k=0}^\infty\frac{||\pmb\xi_{k}||^2 }{\sigma^2_k} }\right\} 
    &\le e^{-\frac{R^2}{2}}.
\end{aligned}
\label{eqn:prob_vio}
\end{equation}
Namely, $\textbf{Rat}$ is bounded by $e^{\epsilon}$ with probability $(1-\delta)$, where $\epsilon = \left\{ \frac{1}{\underline\mu(\mc{L})}(\frac{A}{4} + \frac{R\sqrt{\bar{\mu}(\mc{L})A}}{\sqrt{2}})   \right\}$, $\delta = e^{-\frac{R^2}{2}} $, and $A=\frac{1}{\gamma} \sqrt{\zeta (2(q-p))} ||f_I - f_I'||^2_{\mc{V}_q}$.

Therefore, under proper parameter settings, the mechanism $\mc{M}$ is $(\epsilon, \delta)$-DP with $\epsilon, \delta$ defined above.
\end{proof}

\begin{remark}
\label{rem:1}
Several factors contribute to the privacy parameter $\epsilon$. First, the communication graph topology affects $\bar\mu(\mc{L})$ and $\underline\mu(\mc{L})$. Since $\bar\mu (\mc{L})$ is bounded by $2\triangle(G)$, where $\triangle(G)$ is the maximum degree, and $\underline\mu (\mc{L})$ represents the algebraic connectivity affected by the graph connectivity, a strongly and evenly connected graph helps EFPSN to achieve the best privacy performance. In addition, the parameters $\gamma, q,p$ need to be deliberately chosen to decrease $\epsilon$. Since $\gamma$ is the noise magnitude, it is natural that a larger $\gamma$ contributes to a smaller $\epsilon$ and, consequently,  better privacy.

\end{remark}

\begin{remark}
\label{rem:2}
The term $||f_I - f_I'||^2_{\mc{V}_q}$ measures the privacy-preserving capacity of EFPSN from a different perspective. Given some privacy budget $\epsilon$, the larger it is, the more functions can be protected by our mechanism. Additionally, this term is inversely proportional to $\gamma$, suggesting protecting a larger dataset requires more noise.
\end{remark}

\begin{remark}
\label{rem:3}
Choosing arbitrarily large $R,\gamma$ while keeping $\frac{R}{\sqrt{\gamma}}\sim o(1)$ results in arbitrarily small $\epsilon$ and $\delta$ simultaneously. Namely, we are able to attain a nearly perfectly private mechanism while the algorithmic accuracy is not sacrificed at the same time.
\end{remark}


\section{Numerical Experiments}

In this section, we evaluate the performance of the proposed EFPSN algorithm using numerical experiments. 
We consider both convex and non-convex objective functions, and EFPSN is compared with the non-zero-sum functional perturbation mechanism proposed in \cite{nozari_differentially_2016}. Before displaying the results, we first demonstrate how the orthonormal system $\{e_k\}_{k\in I}$ based on which we generate the coefficients-function mapping $\Phi$, is constructed.

\subsection{Generating the Orthonormal System}
Since the number of elements of an orthonormal basis grows factorially as the number of variables increases, even the simplest real-world applications, such as conducting logistic regression on pictures, require at least hundreds of parameters, making the number of elements in a basis prohibitively large. Consequently, for a function with $M$ variables, we pick $m<M$ variables to perturb. Also, instead of generating an orthonormal basis, we only generate an orthonormal system in $L_2$ with $N$ elements.  

Since the generated orthonormal system must belong to some orthonormal basis, perturbing the orthonormal system is equivalent to perturbing the orthonormal basis with some noise coefficient being $0$. Therefore, all of our previous discussion on accuracy and privacy holds under such a perturbing mechanism. 

The orthonormal system is constructed from the Gram-Schimidt orthonormalization of the Taylor functions. Given the tuple $(K,m,N)$, we randomly generate $N$ elements of Taylor functions of $m$ variables, of which the sum of the order is smaller than or equal to $K$. Then we orthonormalizaing all terms using Gram-Schimidt method.

\begin{figure}[!h]
    \centering
    \includegraphics[width=40mm]{ 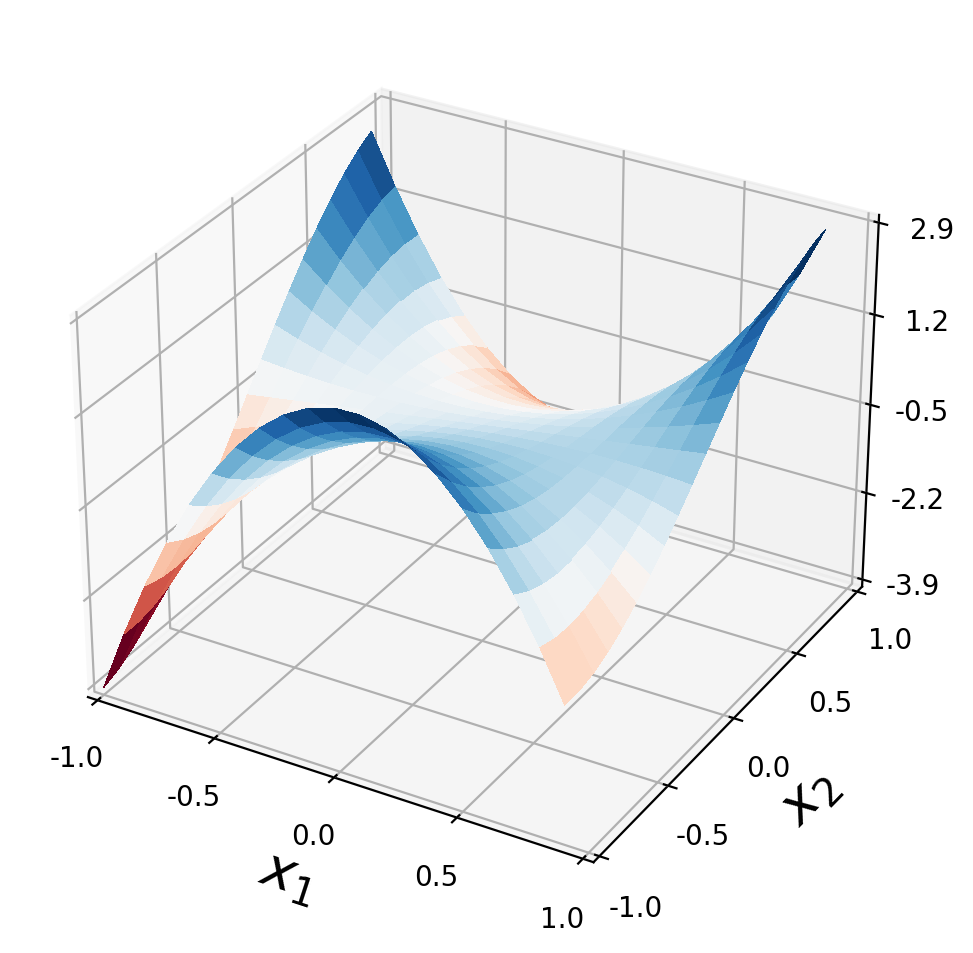}
    \caption{Visualization of a 2d perturbing function.}
    \label{fig:per_func}
\end{figure}

We present an example of one perturbing function when $(K,m,N)=(3,2,5)$ under noise level $\gamma = 1$. The orthonormal system is $\{0.5, 0.866x_2, 3.307 x_2^3-1.984x_2 , 0.866x_1, 2.905x_1^2x_2 - 0.968x_2\}$ and the noise sequence is $\{0.180, 0.628, -0.374, 0.817, 2.015\}$. The corresponding perturbing function is $5.853x_1^2x_2 - 2.137x_2^3 + 0.708x_1 -1.310x_2 + 0.090$. We visualize the perturbing function in Fig. \ref{fig:per_func}.

\subsection{Accuracy Test}

In this part, we validate the accuracy of the proposed EFPSN method.

\subsubsection{Convex Case}
\label{Sec:convex_case}
We consider a classification task on MNIST dataset using logistic regression.

We implement the logistic regression model using PyTorch. Specifically, since the picture in MNIST is of size $28\times 28$ and of channel 1, there is only one linear layer in the model, of which the input and output dimensions are 784 and 10 respectively. Together with bias, the model has 7850 parameters. We set the perturbing parameters $(K,m,N) = (1, 10, 10)$. And the experiments are conducted under different noise level $\gamma\in \{1e^{-2}, 1e^{-1}, 1e^0, 1e^1, 1e^{2}, 1e^{3}, 1e^{4}\}$.

\begin{figure}[!h]
    \centering
    \includegraphics[width=30mm]{ 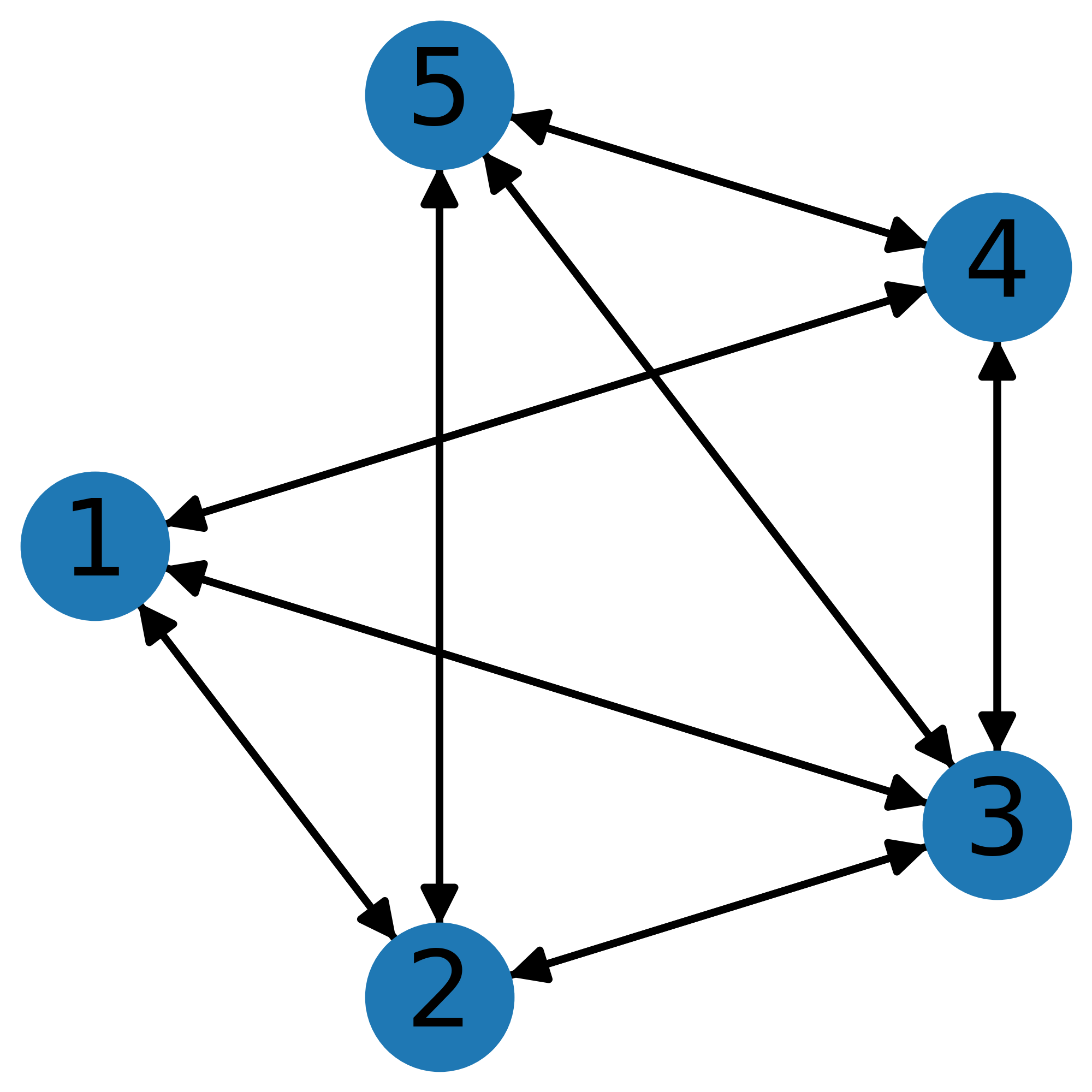}
    \caption{Network structure between agents.}
    \label{fig:network_structure}
\end{figure}

Consider 5 agents in the network, connected as shown in Fig. \ref{fig:network_structure}. Each agent holds the same number of randomly assigned training data points. We adopt decentralized stochastic gradient descent in Phase II in Alg. \ref{alg:EFPSN}. The batch size is set to $64$. The initial learning rate equals $0.2$. Each agent conducts $10000$ gradient updates. The learning rate remains fixed in the first $2000$ steps, and drops to $4e^{-5}$ at the last step.

\begin{figure}
\centering
  \begin{subfigure}[b]{0.35\textwidth}
    \includegraphics[width=\textwidth]{ 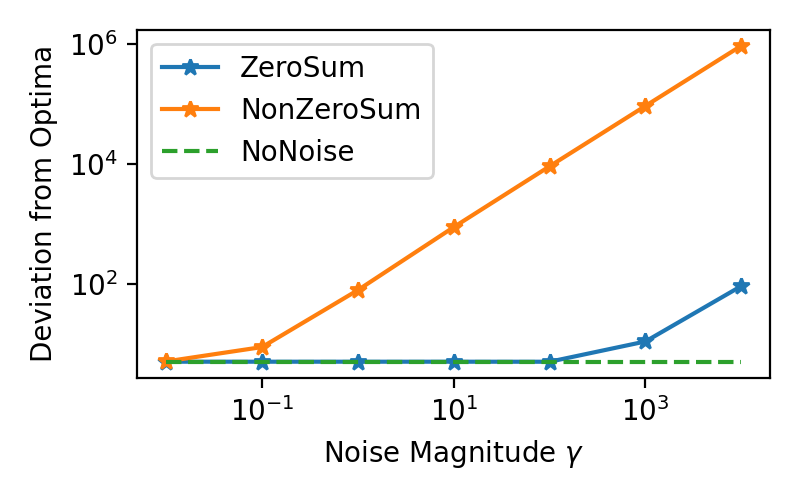}
    \caption{Deviation}
    \label{fig:con_dev}
  \end{subfigure}
 \hfill
  \begin{subfigure}[b]{0.35\textwidth}
    \includegraphics[width=\textwidth]{ 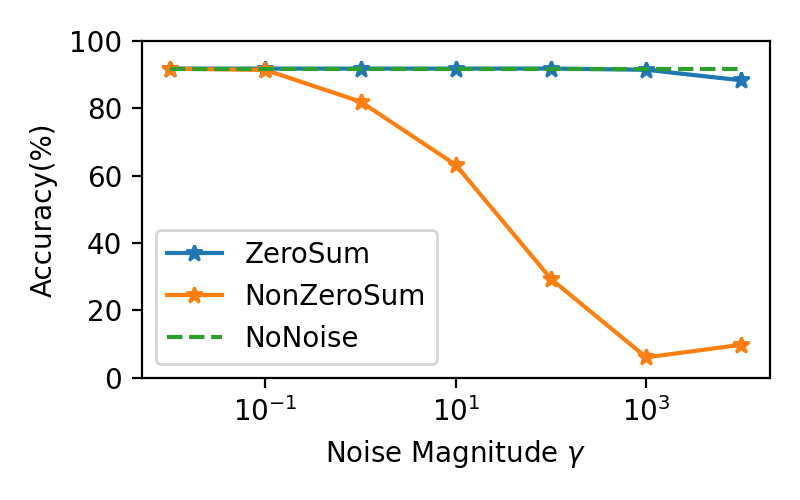}
    \caption{Accuracy}
    \label{fig:con_acc}
  \end{subfigure}
  \caption{Deviation and Accuracy for Logistic Regression}
  \label{fig:con_acc_res}
\end{figure}

The result is shown in Fig. \ref{fig:con_acc_res}. In Fig. \ref{fig:con_dev}, the horizontal axis displays different noise magnitude $\gamma$, and the vertical axis represents the deviation from the optimal solution, which is $||\overline {\mbf{x}} - \mbf{x}^\ast||$. Specifically, $\overline{\mbf{x}}=\frac{1}{n}\sum_{i\in\mc{N}} \mbf{x}_i$, and $\mbf{x}^\ast$ is the solution generated by centralized gradient descent. 

When $\gamma = 1e^{-2}$, the results from EFPSN and the non-zero-sum algorithm are nearly identical, both close to the noise-free case. However, as $\gamma$ increases, the solution obtained from the non-zero-sum method starts to deviate quickly from optima, and such deviation spikes at a roughly constant rate. The blue line (EFPSN solution) does not rise until $\gamma=1e^3$, at which point the orange line is 4 magnitudes higher. 

The rise in the blue line stems from the slight disagreements between local decision variables $\mbf{x}_i$. Note that our perturbed function only guarantees zero-sum when each agent holds the same decision variable. Our previous error analysis assumes such exact agreement between agents. Yet, in practice, slight differences between $\mbf{x}_i$ exist due to finite step size. With EFPSN, one can always generate a more accurate solution using a finer learning rate.

Fig. \ref{fig:con_acc} demonstrates the classification accuracy of the logistic model under different algorithms and noise levels. The pattern matches that of Fig. \ref{fig:con_dev}. Basically, for the non-zero-sum algorithm, the test accuracy starts to drop dramatically when $\gamma$ reaches $1e^0$. In contrast, the model trained by EFPSN remains as accurate as the noise-free case until $\gamma = 1e^4$. Namely, EFPSN provides much more (at least 4 magnitudes larger) privacy budgets while not degenerating accuracy.

\subsubsection{Non-convex Case}
To further justify our method under the non-convex setting, we consider an image classification task on MNIST using Convolutional Neural Network. We adopt the classic LeNet \cite{lecun_gradient-based_1998}, which consists of 13426 parameters. We choose to perturb the bias of its last linear layer, and again we set $(K,m,N)=(1,10,10)$. Except for the model, all the settings are identical to the convex case.

\begin{figure}
\centering
  \begin{subfigure}[b]{0.35\textwidth}
    \includegraphics[width=\textwidth]{ 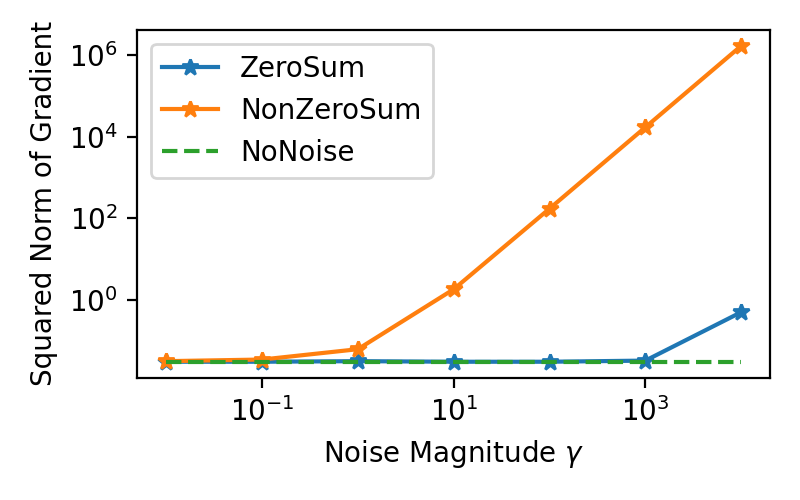}
    \caption{Squared Norm of Average Gradient}
    \label{fig:noncon_dev}
  \end{subfigure}
 \hfill
  \begin{subfigure}[b]{0.35\textwidth}
    \includegraphics[width=\textwidth]{ 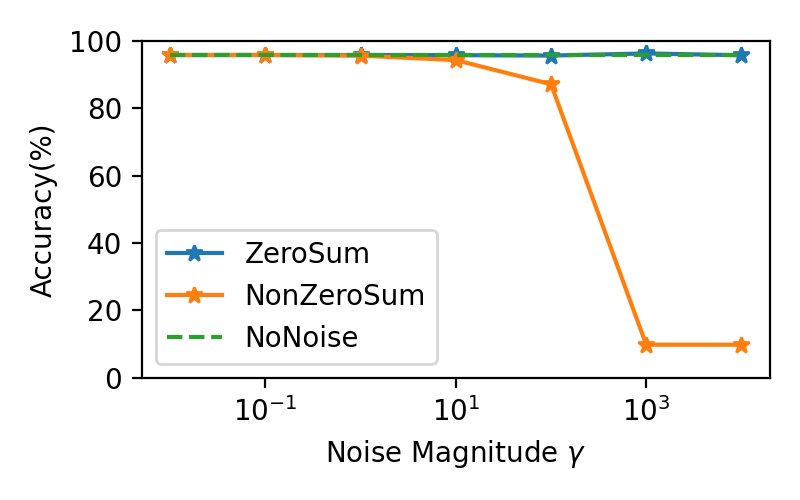}
    \caption{Accuracy}
    \label{fig:noncon_acc}
  \end{subfigure}
  \caption{Squared Norm of Average Gradient and Accuracy for LeNet}
  \label{fig:noncon_acc_res}
\end{figure}

In fig. \ref{fig:noncon_dev}, the y-axis depicts the squared norm of the average gradient of all the agents, $||\frac{\sum_{i\in\mc{N}} \nabla f_i(x_i^k)}{|\mc{N}|}||^2$. Trends similar to the convex case reappear in the non-convex case. Unlike the non-zero-sum method, EFPSN generates solutions closer to the stationary point under all noise levels. Also, in fig. \ref{fig:noncon_acc}, EFPSN remains as accurate as the noise-free case even when $\gamma=1e^4$. The non-zero-sum method, nonetheless, could barely hold its accuracy as $\gamma=1e^0$. 

The experiments imply the efficacy of our method under both convex and non-convex problems.

\subsection{Privacy Test}

To further validate EFPSN's efficacy in privacy-preserving, we conduct  DLG \cite{zhu_deep_nodate} attack on agents with zero-sum and non-zero-sum noise, respectively. The general idea of DLG is to construct some dummy data and try to match its gradient with the ground truth. The algorithm is shown in alg. \ref{alg:DLG}. For better results, we use iDLG \cite{DBLP:journals/corr/abs-2001-02610}, a more efficient and stable version of DLG.

\begin{algorithm}
\caption{Deep Leakage from Gradients}\label{alg:DLG}
\begin{algorithmic}[1]
\Require $F(\mbf{x},W)$: Differentiable model; $W$: parameter weights; $\nabla W$ : gradients calculated by training data
\Ensure Private training data $\mbf{x}, \mbf{y}$

\State{$\mbf{x}'_1 \leftarrow N(0,1), \mbf{y}'_1 \leftarrow N(0,1)$}

\For {$i\leftarrow 1 $ to $n$}
    
    \State{$\nabla W_i' \leftarrow \partial l (F(\mbf{x}'_i,W_t), \mbf{y}_i')/\partial W_t$}
    \State{$\mbb{D}_i \leftarrow ||\nabla W_i' - \nabla W||^2 $}
    \State{$\mbf{x}_{i+1}' \leftarrow  \mbf{x}_{i}' - \alpha \nabla_{\mbf{x}_{i+1}'} \mbb{D}_i$}
    \State{$\mbf{y}_{i+1}' \leftarrow  \mbf{y}_{i}' - \alpha \nabla_{\mbf{y}_{i+1}'} \mbb{D}_i$}

\EndFor
\State \Return{$\mbf{x}_{n+1}', \mbf{y}_{n+1}'$}
\end{algorithmic}
\end{algorithm}

With either EFPSN or the non-zero-sum noise method, each agent receives some non-zero functional perturbation of roughly the same magnitude. Since iDLG is carried out at the agent level, it makes no difference between the two algorithms. Therefore, we conduct iDLG on agent 1 in fig. \ref{fig:network_structure}, with the noise generated from EFPSN at different noise level ($\gamma \in \{1e^1, 1e^2, 1e^3, 1e^4\}$).

Specifically, we assume the mixing matrix is known to the attacker. And the attacker has access to at least one of the communication channels connected to agent 1 (either eavesdropping or corrupting 1's neighbor will do). Therefore, the attacker knows agent 1's perturbed gradient $\nabla \hat f_1(x_1^k)$ and agent 1's decision parameters $x_1^k$. The attacker does not know the functional perturbation $\Phi(\bar\eta_i)$. Consequently, the true gradient  $\nabla f_1(x_1^k)$ remains unrevealed. Namely, the attacker is trying to recover the raw data using inexact gradient information. And the larger the noise level $\gamma$, the more inexact the gradient is. 

\begin{figure}
\centering
  \begin{subfigure}[b]{0.10\textwidth}
    \includegraphics[width=\textwidth]{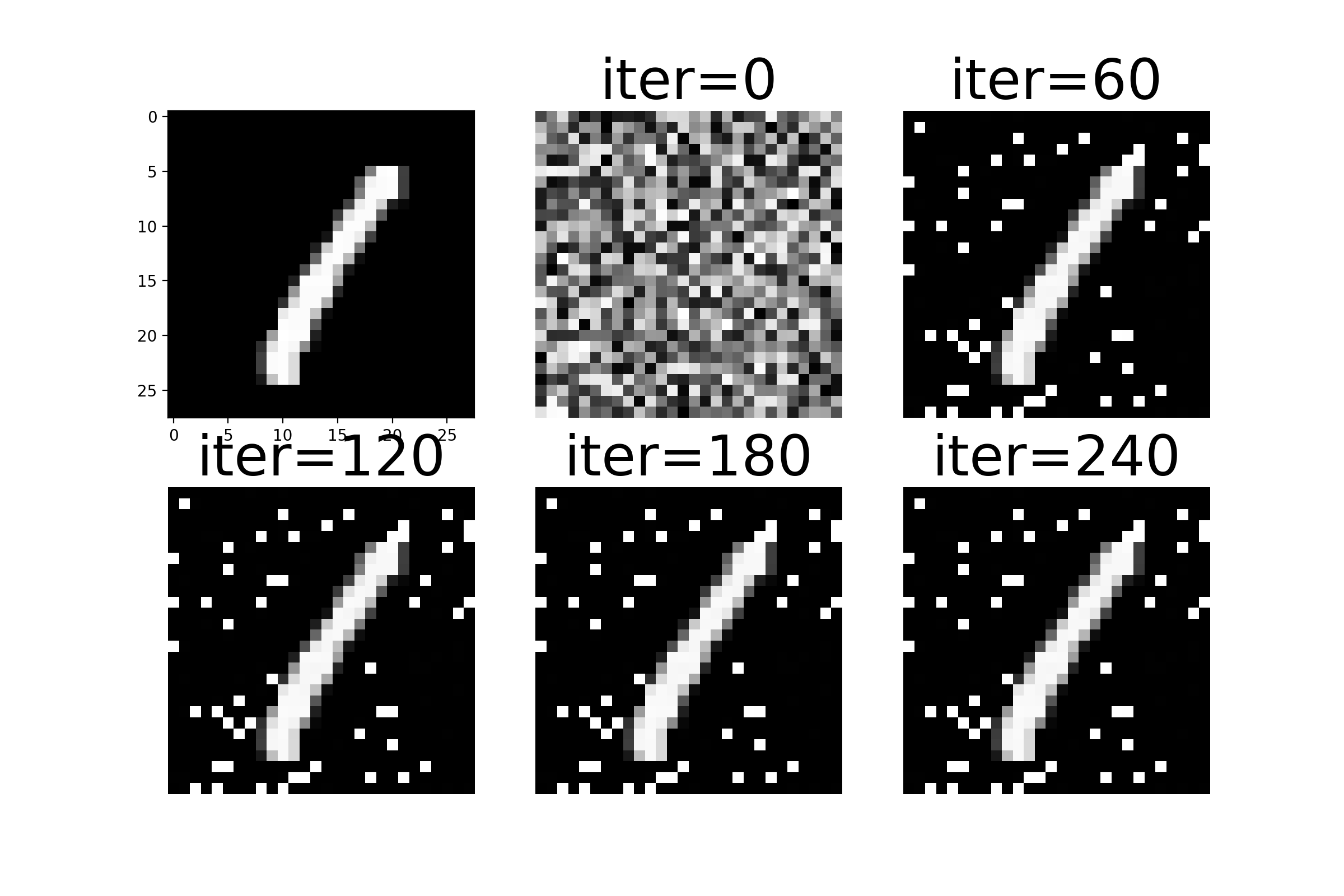}
    \caption{$\gamma=1e^1$}
    \label{fig:iDLG_log_1}
  \end{subfigure}
  \begin{subfigure}[b]{0.10\textwidth}
    \includegraphics[width=\textwidth]{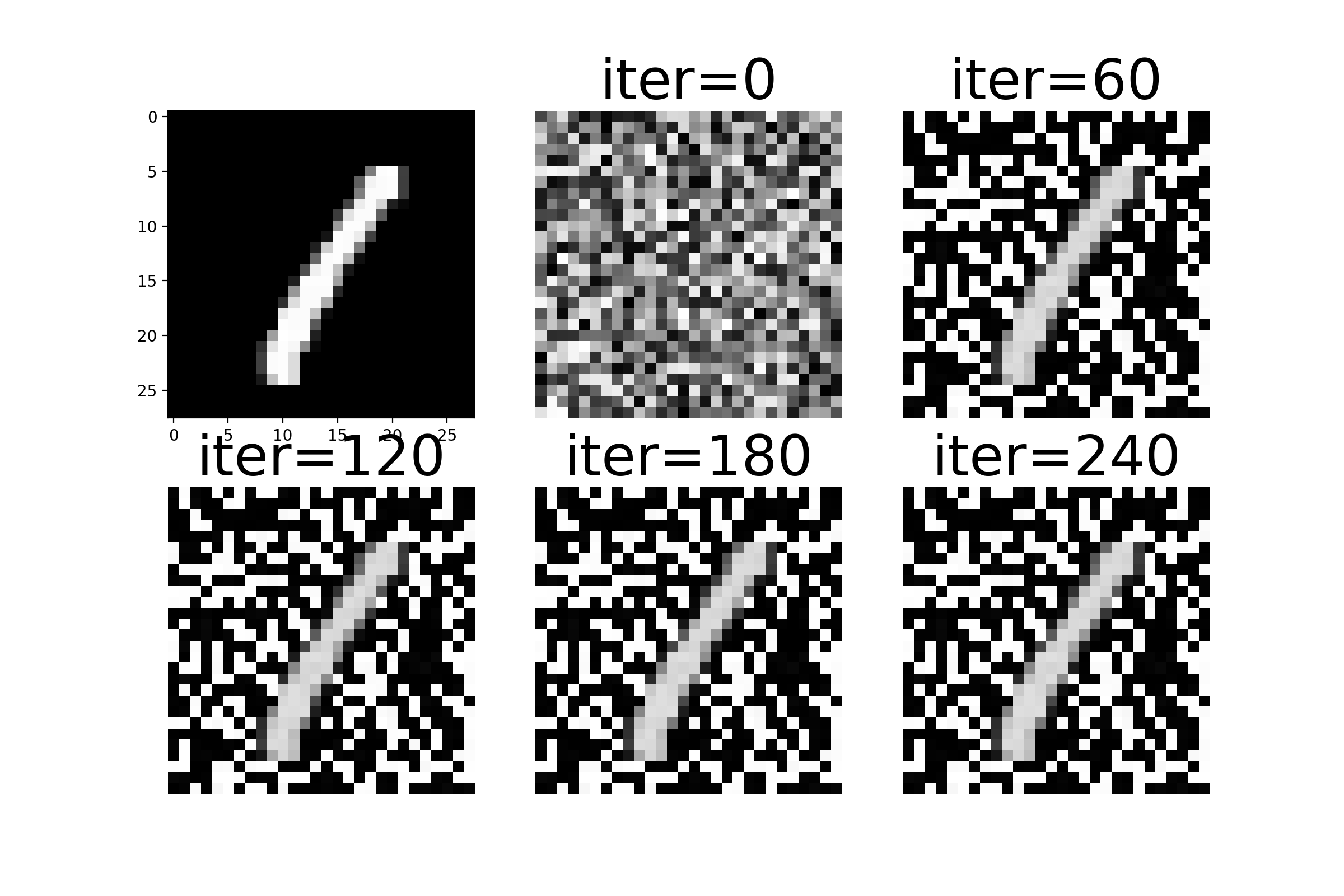}
    \caption{$\gamma=1e^2$}
    \label{fig:iDLG_log_2}
  \end{subfigure}
  \begin{subfigure}[b]{0.10\textwidth}
    \includegraphics[width=\textwidth]{ 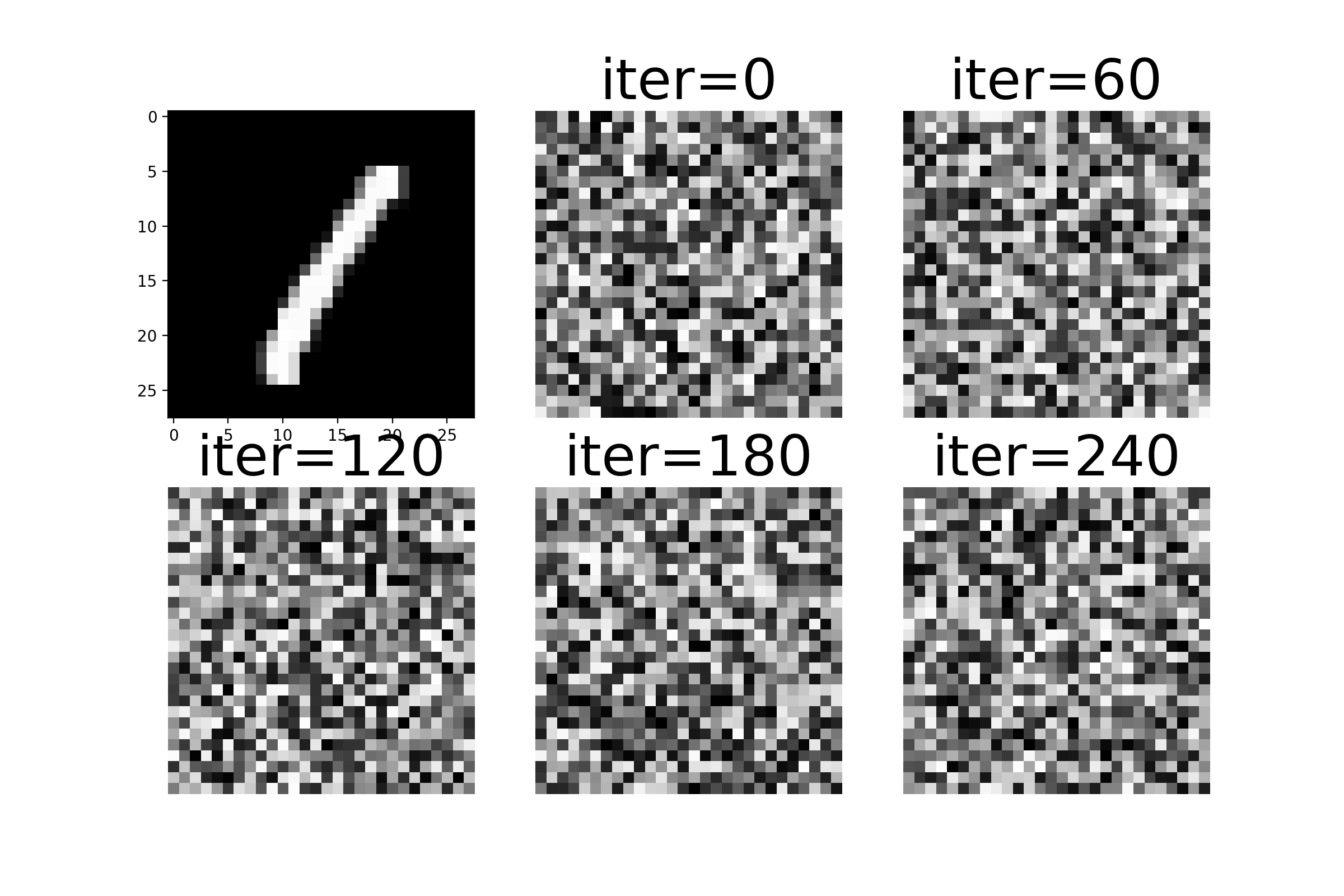}
    \caption{$\gamma=1e^3$}
    \label{fig:iDLG_log_3}
  \end{subfigure}
  \begin{subfigure}[b]{0.10
  \textwidth}
    \includegraphics[width=\textwidth]{ 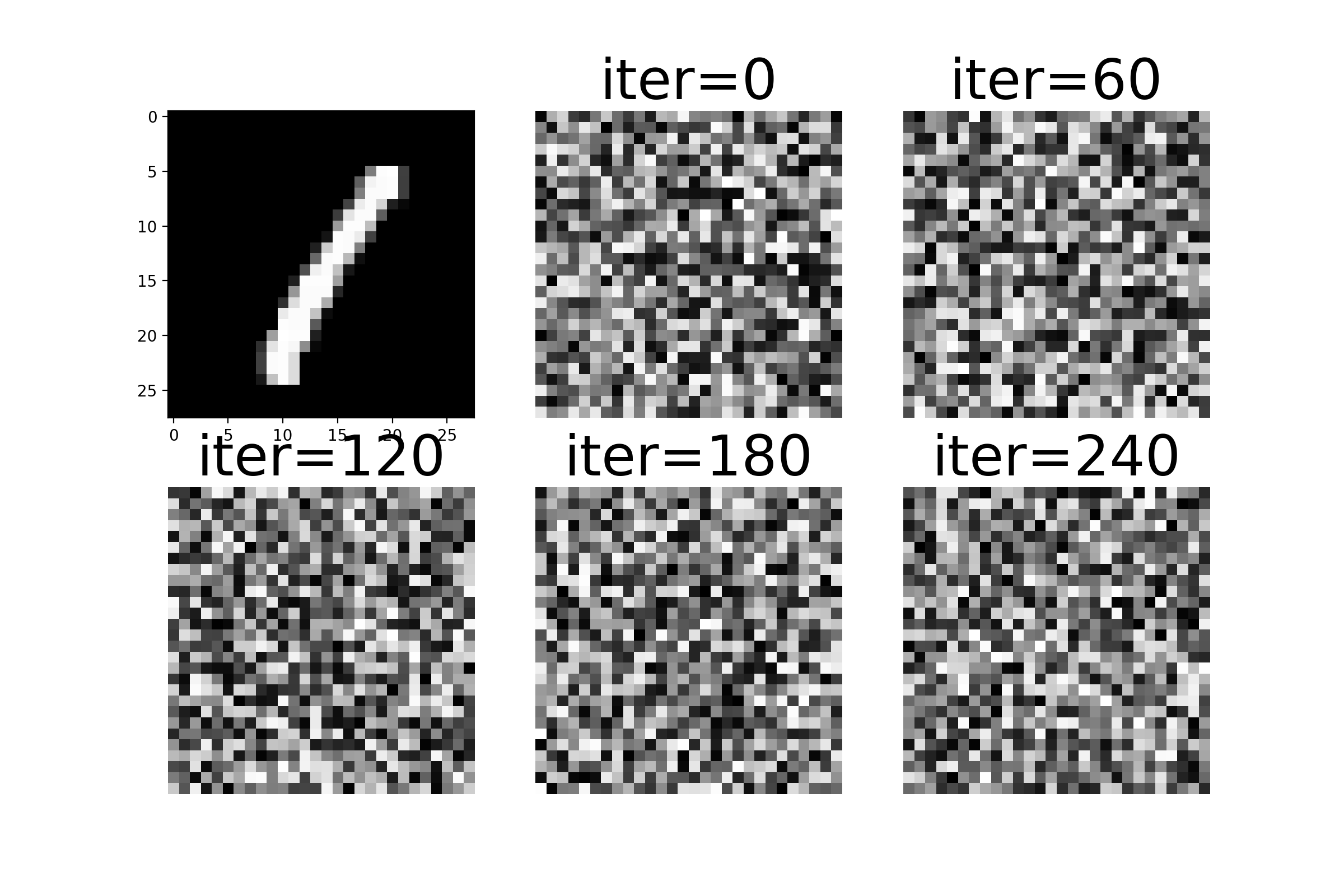}
    \caption{$\gamma=1e^4$}
    \label{fig:iDLG_log_4}
  \end{subfigure}
  \caption{iLDG attacker's inference results on Logistic Regression}
  \label{fig:iDLG_log}
\end{figure}

\begin{figure}
\centering
  \begin{subfigure}[b]{0.10\textwidth}
    \includegraphics[width=\textwidth]{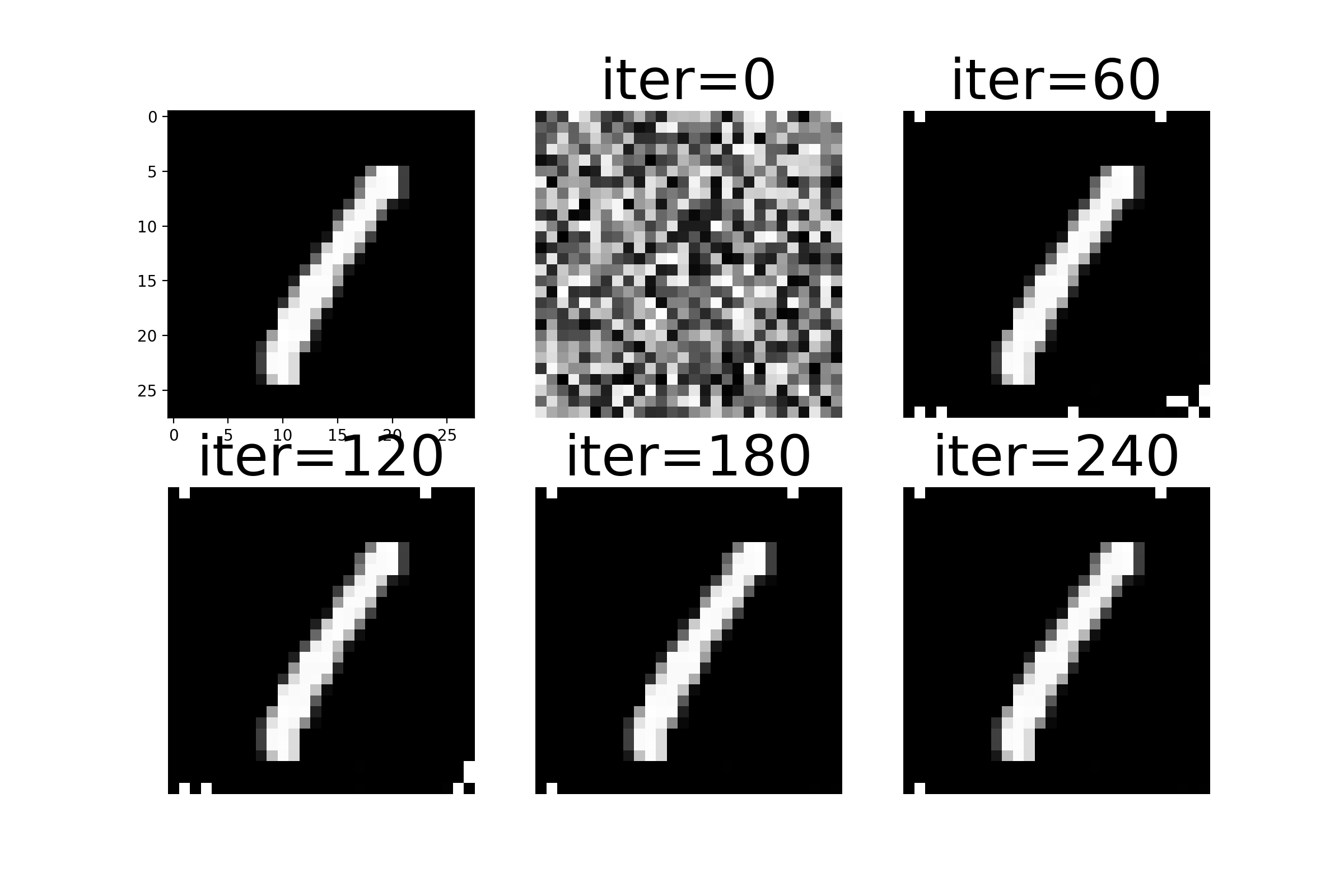}
    \caption{$\gamma=1e^1$}
    \label{fig:iDLG_lenet_1}
  \end{subfigure}
  \begin{subfigure}[b]{0.10\textwidth}
    \includegraphics[width=\textwidth]{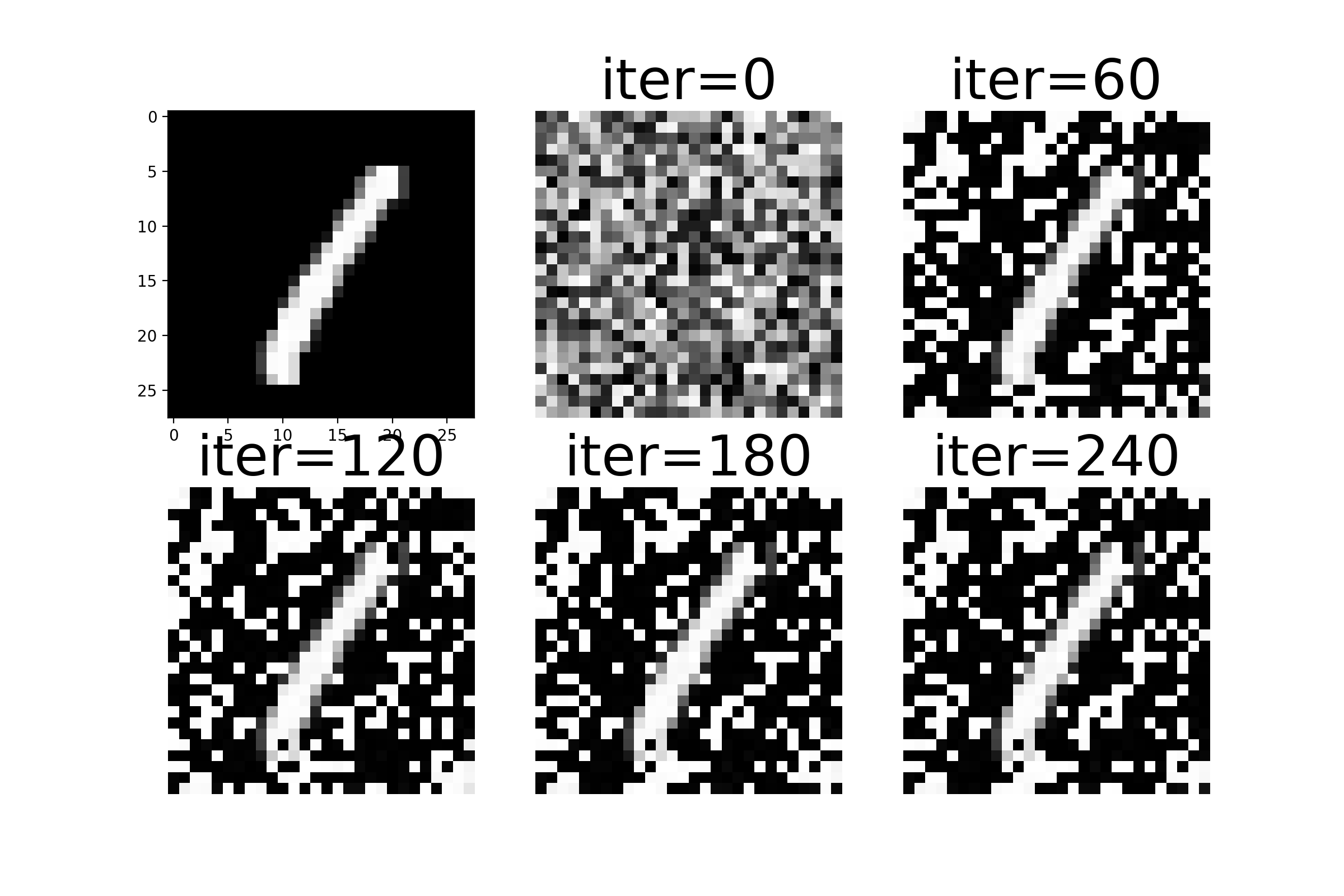}
    \caption{$\gamma=1e^2$}
    \label{fig:iDLG_lenet_2}
  \end{subfigure}
  \begin{subfigure}[b]{0.10\textwidth}
    \includegraphics[width=\textwidth]{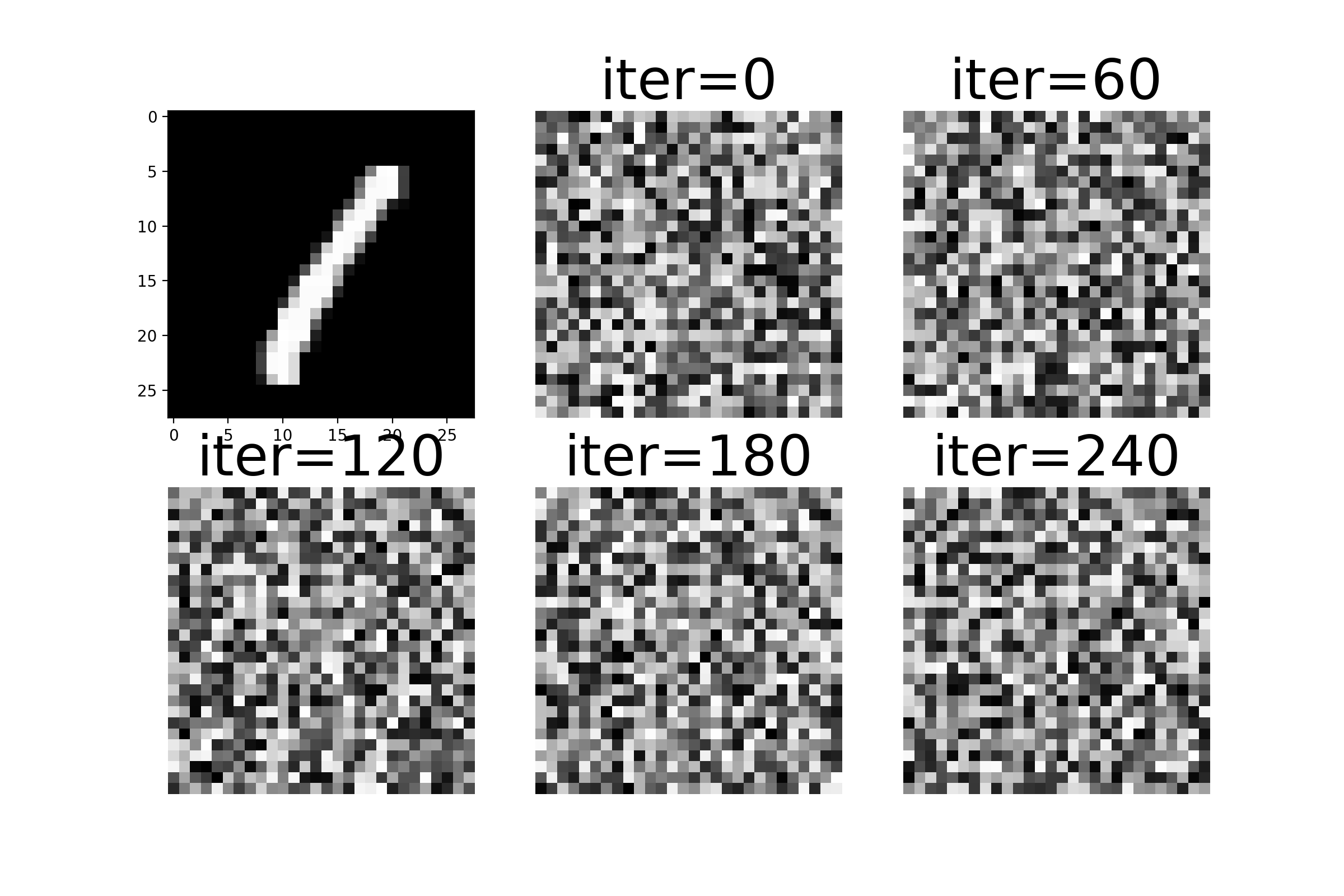}
    \caption{$\gamma=1e^3$}
    \label{fig:iDLG_lenet_3}
  \end{subfigure}
  \begin{subfigure}[b]{0.10
  \textwidth}
    \includegraphics[width=\textwidth]{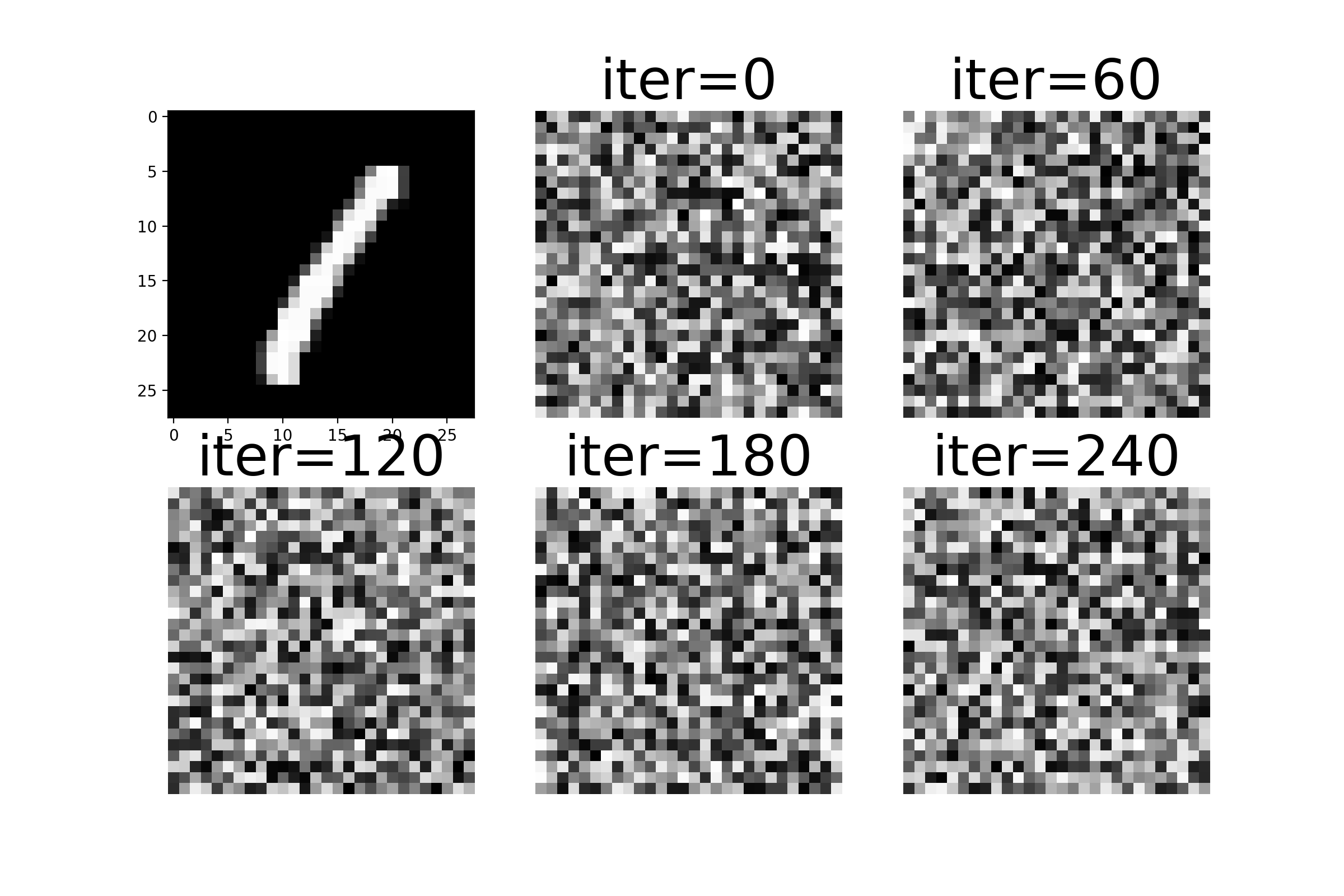}
    \caption{$\gamma=1e^4$}
    \label{fig:iDLG_lenet_4}
  \end{subfigure}
  \caption{iLDG attacker's inference results on LeNet}
  \label{fig:iDLG_lenet}
\end{figure}

Fig. \ref{fig:iDLG_log} and fig. \ref{fig:iDLG_lenet} represent the iDLG attacker's typical inference result on the Logistic model and LeNet. The top left subfigure is the raw data. And the remaining subfigures are the adversary's estimate of the raw data at different iterations (from 0 to 240). As $\gamma$ increases, the retrieved picture becomes blurred. Interestingly, though we are perturbing the original problem functionally, it is equivalent to directly perturbing the dataset. Generally, after $\gamma \ge 1e^3$, the recovered picture is unrecognizable for humans. When $\gamma \ge 1e^3$, however, the accuracy of the model trained by the non-zero-sum method has dropped below 10\% for both convex and non-convex problems (as shown in fig. \ref{fig:con_acc_res} and fig. \ref{fig:noncon_acc_res}). This suggests that the non-zero-sum solution would be too inaccurate to provide enough privacy. In contrast, the EFPSN solution has comparable accuracy to the noise-free case. Thus, EFPSN is capable of preserving privacy without degenerating accuracy.

\section{Conclusion}

In the paper, we proposed the Encrypted Functional Perturbation with Structured Noise algorithm that solves the decentralized optimization problem \ref{eqn:problem_definition} privately and accurately. Given exact consensus between agents, EFPSN could eliminate the privacy-accuracy trade-off by constructing a zero-sum functional perturbation. Since such construction requires secure communication between agents, we adopt the Paillier encryption scheme to fight against eavesdropping attackers. We rigorously proved the privacy property of EFPSN under the differential privacy framework. Simulations confirmed the efficacy of EFPSN in protecting privacy while maintaining accuracy.

\bibliographystyle{ieeetr}
\bibliography{ref.bib}
\end{document}